\documentclass[a4paper, 11pt]{article}

\usepackage[T1]{fontenc}
\usepackage[utf8]{inputenc}
\usepackage{lmodern}
\usepackage{xcolor}
\usepackage{amsmath}
\numberwithin{equation}{section}
\usepackage{amsthm}
\usepackage{amssymb}
\usepackage{physics}
\usepackage{mathtools}
\usepackage{mathrsfs}
\usepackage[backend=bibtex,style=alphabetic]{biblatex}
\usepackage[linkcolor=black, urlcolor=black, citecolor=blue, colorlinks=true]{hyperref}

\usepackage[
  a4paper,
  left=2.8cm,
  right=2.8cm,
  top=3cm,
  bottom=3cm
]{geometry}

\newcommand{\C}{\mathbb{C}}

\newcommand{\N}{\mathbb{N}}

\newcommand{\R}{\mathbb{R}}

\newcommand{\1}{\mathbf{1}}

\newcommand{\cB}{\mathcal{B}}

\newcommand{\cU}{\mathcal{U}}
\newcommand{\cL}{\mathcal{L}}

\newcommand{\cT}{\mathcal{T}}

\newcommand{\mi}{\mathrm{i}}
\newcommand{\me}{\mathrm{e}}
\newcommand{\md}{\mathrm{d}}

\newtheorem{theorem}{Theorem}[section]
\newtheorem{corollary}[theorem]{Corollary}
\newtheorem{proposition}[theorem]{Proposition}
\newtheorem{lemma}[theorem]{Lemma}

\newtheorem*{assumption*}{Assumption}

\theoremstyle{definition}

\newtheorem{example}[theorem]{Example}
\newtheorem{remark}[theorem]{Remark}

\DeclareMathOperator{\dom}{Dom}
\DeclareMathOperator{\gap}{gap}
\DeclareMathOperator{\myspan}{span}
\DeclareMathOperator{\diag}{diag}

\providecommand{\keywords}[1]{\par\smallskip\noindent{\small \textbf{Keywords:} #1}}
\providecommand{\subjclass}[1]{\par\smallskip\noindent{\small \textbf{2020 Mathematics Subject Classification.} #1}}

\title{On the Existence of the KMS Spectral Gap in Gaussian Quantum Markov Semigroups}

\author{
Zheng Li
\\[1ex]
\small
School of Mathematics and Statistics, Central South University,
Changsha, 410083, China.
\\[0.5ex]
\small
\href{mailto:zheng.li@csu.edu.cn}{zheng.li@csu.edu.cn}
}

\date{}

\allowdisplaybreaks

\begin{document}

    \maketitle

    \begin{abstract}
        In \cite{fagnola2025spectral}, it was shown that a Gaussian quantum Markov semigroup on the $d$-mode bosonic Fock space with a unique faithful normal invariant state has a positive GNS spectral gap if and only if the matrix $[U\,\overline V]$, formed from the coefficients of the Kraus operators, has $2d$ linearly independent columns. In this paper, we establish the corresponding criterion for the KMS spectral gap. After standardizing the invariant Gaussian state and grouping the modes according to their inverse temperatures, let $U_n$ and $V_n$ denote the submatrices of $U$ and $V$, respectively, corresponding to the $n$-th equal-temperature block. We prove that $\ker U_n=\ker V_n$ for every $n$, and that the KMS spectral gap is positive if and only if this common kernel is trivial for every equal-temperature block. Thus, whereas the GNS spectral gap requires a global maximal-rank condition on the Kraus coefficients, the KMS spectral gap requires maximal column rank only within each equal-temperature block.
    \end{abstract}



    \section{Introduction}

    A quantum Markov semigroup (QMS) on a von Neumann algebra is a one-parameter semigroup of normal, completely positive, identity-preserving maps. Its generator is the noncommutative analogue of a Markov generator \cite{gorini1976completely,lindblad1976generators}. In the present paper, the underlying algebra is $\cB(\mathsf h)$, where $\mathsf h$ is a finite-mode bosonic Fock space.

    A QMS is called Gaussian when its predual semigroup preserves Gaussian normal states \cite{caruso2008multi,heinosaari2010semigroup,poletti2022characterization}. Such semigroups are noncommutative analogues of Ornstein--Uhlenbeck semigroups. Their action on Weyl operators and the evolution of Gaussian covariances are governed by finite-dimensional drift and diffusion operators; see \cite{poletti2022characterization,agredo2022decoherence,fagnola2025spectral}. Consequently, many analytic questions reduce to matrix identities and continuous Lyapunov equations. General references on bosonic Gaussian states and maps include \cite{holevo2011probabilistic,weedbrook2012gaussian,serafini2023quantum}.

    Let $(\cT_t)_{t\geq 0}$ be a Gaussian QMS on $\cB(\mathsf h)$, where $\mathsf h=\Gamma(\C^d)$ is the $d$-mode bosonic Fock space. We assume that $(\cT_t)_{t\geq 0}$ admits a unique faithful normal invariant state $\rho$. Its generator $\cL$ has the generalized Gorini--Kossakowski--Sudarshan--Lindblad (GKSL) form \cite{poletti2022characterization}

    \begin{equation}\label{eq-gksl-generator}
        \cL(x) = \mi[H,x] - \frac{1}{2}\sum_{\ell=1}^m
        \left(L_\ell^*L_\ell x-2L_\ell^*xL_\ell+xL_\ell^*L_\ell\right),
        \qquad x\in\dom\cL,
    \end{equation}
    where
    \begin{align}
        H &= \sum_{j,k=1}^{d} \left( \Omega_{jk} a_j^\dagger a_k + \frac{\kappa_{j k}}{2} a^\dagger_j a^\dagger_k + \frac{\overline{\kappa_{jk}}}{2} a_j a_k \right) + \sum_{j=1}^d \left( \frac{\zeta_j}{2} a^\dagger_j + \frac{\overline{\zeta_j}}{2} a_j \right), \label{eq-Hamiltonian} \\ 
        L_\ell &= \sum_{k=1}^d \left( \overline{v_{\ell k}} a_k + u_{\ell k} a^\dagger_k \right), \label{eq-Kraus-operator}
    \end{align}
    Here $a_j$ and $a_j^\dagger$ are the annihilation and creation operators on the $j$-th mode, respectively; their precise definitions are given in \eqref{eq-definition-annihilation-operator} and \eqref{eq-definition-creation-operator}. We call $H$ the \emph{Hamiltonian} and $L_1,\ldots,L_m$ the \emph{Kraus operators}. Throughout the paper, we fix a minimal Kraus family. After removing linear redundancies, one may therefore assume that $1\leq m\leq 2d$ (see \cite[Definition 2.1]{agredo2022decoherence}). This is a choice of representation rather than an additional hypothesis on the Gaussian QMS. The brackets $[\,\cdot\,,\,\cdot\,]$ and $\{\,\cdot\,,\,\cdot\,\}$ denote the commutator and anticommutator, respectively.

    We collect the quadratic coefficients of the Hamiltonian in the matrices
    \begin{equation} \label{eq-definition-omega-kappa}
        \Omega := \left[ \Omega_{j k} \right]_{1 \le j, k \le d} = \Omega^*, \quad \kappa := \left[ \kappa_{j k} \right]_{1 \le j, k \le d} = \kappa^T,
    \end{equation}
    where $\cdot^{T}$ denotes transpose and $\cdot^*$ conjugate transpose. Thus, $\Omega$ is Hermitian and $\kappa$ is symmetric. We collect the coefficients of the Kraus operators in the matrices
    \begin{equation} \label{eq-definition-U-and-V}
        V := \left[ v_{\ell k} \right]_{1 \le \ell \le m , 1 \le k \le d }, \quad U := \left[ u_{\ell k} \right]_{1 \le \ell \le m,  1 \le k \le d},
    \end{equation}
    and $\zeta=(\zeta_j)_{j=1}^d\in\C^d$. Thus, $U,V\in M_{m\times d}(\C)$. Square brackets around a single array denote a matrix and will not be confused with the two-variable commutator notation.

    Since $\rho$ need not be tracial, different embeddings of $\cB(\mathsf h)$ into the Hilbert--Schmidt space $\cB_2(\mathsf h)$ lead to different induced semigroups. For $s\in[0,1]$, the induced contraction semigroup $(T_t^{(s)})_{t\geq0}$ is determined on the dense subspace $\rho^{s/2}\cB(\mathsf h)\rho^{(1-s)/2}$ by
    \begin{equation} \label{eq-induced-semigroup-and-qms} 
        T_t^{(s)}\bigl(\rho^{s/2}X\rho^{(1-s)/2}\bigr)
        :=\rho^{s/2}\cT_t(X)\rho^{(1-s)/2},
        \qquad X\in\cB(\mathsf h).
    \end{equation} 
    The cases $s=0$, $s=1/2$, and $s=1$ are called the GNS, KMS, and anti-GNS embeddings, respectively. We denote the generator of $(T_t^{(s)})_{t\geq0}$ by $L^{(s)}$. Our main object is $L^{(1/2)}$.

    The choice of embedding affects both symmetry and coercivity properties; see, for example, \cite{carlen2017gradient,wirth2024christensen,vernooij2023derivations}. We use the following numerical-range definition, which is consistent with the classical spectral-gap formulation in \cite[Theorem 2.3]{liggett1989exponential}. This definition is also equivalent to that of a spectral gap in \cite[Definition 1]{fagnola2025spectral}; see \cite[Appendix A]{GIROTTI2026130150} for a proof of this equivalence. If $L$ is a maximally dissipative operator on a Hilbert space, set

    \begin{equation} \label{eq-definition-spectral-gap}
        \gap L := \inf \left\{ -\Re\langle x,Lx\rangle :
        x\in\dom L\cap(\ker L)^\perp,\ \norm{x}=1 \right\}.
    \end{equation}
    We say that $L$ has a \emph{spectral gap} if $\gap L>0$. This quantity is a coercivity constant for the generator on the orthogonal complement of its kernel. In appropriate invariant subspaces it yields exponential decay of the induced semigroup and is closely related to Poincar\'e-type inequalities. For related results on entropy decay and logarithmic Sobolev inequalities, see \cite{carbone2008hypercontractivity,carbone2015logarithmic,bardet2022hypercontractivity,bertini2024trace}.
    
    The spectral gap of symmetric one-mode Gaussian generators was studied in \cite{cipriani2000spectral}. For multimode generators without a symmetry assumption, \cite{fagnola2025spectral} proves that the GNS spectral gap is positive exactly when the $m\times 2d$ matrix $[U\,\overline V]$ has full column rank. That work also gives an example with a positive KMS spectral gap and a vanishing GNS spectral gap. In one mode, positivity of the KMS spectral gap follows under the present invariant-state assumption \cite{GIROTTI2026130150,fagnola_li_2025spectral}.

    \medskip
    \noindent\textbf{Main result.}
    After standardizing the QMS as in Proposition \ref{proposition-parameters-standardized-GQMS}, we group the modes according to the inverse temperatures in \eqref{eq-definition-inverse-temperatures}. For each distinct inverse temperature $\beta_{(n)}$, the set $I_n$ in \eqref{eq-definition-I-n} is the corresponding equal-temperature block. The matrices $U_n$ and $V_n$ defined in \eqref{eq-definition-U-V-n} are obtained by retaining the columns of $U$ and $V$, respectively, whose indices belong to $I_n$. Lemma \ref{lemma-same-block-relation-U-and-V} shows that $\ker U_n=\ker V_n$, and Theorem \ref{theorem-main-2} proves that the KMS spectral gap is positive if and only if this common kernel is trivial for every $n$. Equivalently, $U_n$ and $V_n$ have full column rank on every equal-temperature block. Thus, the KMS criterion is blockwise, whereas the GNS criterion requires the global matrix $[U\,\overline V]$ to have full column rank. As an immediate consequence, a positive GNS spectral gap implies a positive KMS spectral gap; see Corollary \ref{corollary-GNS-implies-KMS}. Moreover, unlike in the one-mode case, a multimode Gaussian QMS satisfying the same invariant-state assumption need not have a positive KMS spectral gap; see Example \ref{example-two-modes-the-same-inverse-temperatures}.

    \medskip

    The paper is organized as follows. Section \ref{section-introduction-of-GQMS} collects the required facts about Gaussian states and Gaussian QMSs. Section \ref{section-existence-kms-spectral-gap} proves the blockwise criterion and compares it with the GNS condition. Section \ref{section-application-boson-chain-model} applies the criterion to a three-mode bosonic chain. Appendix \ref{section-appendix-matrix-identification} fixes our conventions for real-linear operators. Appendix \ref{section-appendix-invariance-spectral-gap} proves invariance under unitary changes of the minimal Kraus family and under the choice of Williamson diagonalization.

    \section{Preliminaries on Gaussian Quantum Markov Semigroups} \label{section-introduction-of-GQMS}

    We begin with the notation and structural results used in the proof of the main theorem.

    Let $\mathsf h:=\Gamma(\C^d)$ be the $d$-mode bosonic Fock space. Under the canonical identification
    $\Gamma(\C^d)\simeq\Gamma(\C)^{\otimes d}$, each factor $\Gamma(\C)$ is isomorphic to $\ell^2(\N)$. We denote the canonical orthonormal basis of $\mathsf h$ by
    \begin{equation*}
        \left\{ e(n_1,\ldots,n_d):=e_{n_1}\otimes\cdots\otimes e_{n_d}:
        n_1,\ldots,n_d\in\N \right\},
    \end{equation*}
    where $(e_n)_{n\in\N}$ is the canonical basis of $\Gamma(\C)$.
    
    The creation, annihilation, and number operators on the $j$-th mode are denoted by $a_j^\dagger$, $a_j$, and $N_j$, respectively. On the canonical basis,
    \begin{align}
        a_j \, e(n_1, \cdots, n_d) &= \sqrt{n_j} \, e(n_1, \cdots, n_{j-1}, n_{j} - 1, n_{j+1}, \cdots, n_d), \label{eq-definition-annihilation-operator} \\
        a_j^\dagger \, e(n_1, \cdots, n_d) &= \sqrt{n_j + 1} \, e(n_1, \cdots, n_{j-1}, n_{j} + 1, n_{j+1}, \cdots, n_d), \label{eq-definition-creation-operator} \\ 
        N_j e(n_1, \cdots, n_d) &= n_j \, e(n_1, \cdots, n_d). \label{eq-number-operator}
    \end{align}
    Here the first expression is understood to be zero when $n_j=0$.
    
    We also use the momentum and position operators
    \begin{equation*}
        p_j := \mi \, \frac{a^\dagger_j - a_j}{\sqrt{2}}, \quad q_j := \frac{a^\dagger_j + a_j}{\sqrt{2}}.
    \end{equation*}
    For $z\in\C^d$, set
    \begin{equation*}
        a (z) := \sum_{j=1}^d \overline{z_j} a_j, \quad a^\dagger (z) := \sum_{j=1}^d z_j a_j^\dagger,
    \end{equation*}
    and
    \begin{equation} \label{eq-momentum-operator}
        p (z) := \mi \, \frac{a^\dagger (z) - a(z)}{\sqrt{2}} = \mi \, \frac{ \sum_{j=1}^d ( z_j a_j^\dagger - \overline{z_j} a_j ) }{\sqrt{2}}.
    \end{equation}

    The Weyl operator associated with $z\in\C^d$ is
    \begin{equation*}
        W(z) := \me^{ - \mi \sqrt{2} \, p(z) }.
    \end{equation*}
    Our convention gives the Weyl relation
    \begin{equation} \label{eq-weyl-operator-commutation-relation}
        W (z) W(w) = \me^{-\mi \Im \langle z, w \rangle} W(z + w), \quad \forall z, w \in \C^d.
    \end{equation}

    We write $\cB(\mathsf h)$, $\cB_1(\mathsf h)$, and $\cB_2(\mathsf h)$ for the bounded, trace-class, and Hilbert--Schmidt operators on $\mathsf h$, respectively. The Hilbert--Schmidt inner product is denoted by $\langle\cdot,\cdot\rangle_2$. Every normal state on $\cB(\mathsf h)$ is represented by a unique positive trace-class operator of trace one. We use the same symbol for the state and its density operator.

    The \emph{characteristic function} of a normal state $\rho$ is
    \begin{equation*}
        z \mapsto \Tr (\rho W(z)), \quad z \in \C^d.
    \end{equation*}
    A normal state is uniquely determined by its characteristic function \cite[Proposition 2.4]{parthasarathy2010gaussian}; see also \cite[Theorem 5.3.3]{holevo2011probabilistic}.

    A normal state $\rho$ is \emph{Gaussian} if there exist a vector $\omega\in\C^d$ and a positive self-adjoint real-linear operator $S$ on the real Hilbert space $(\C^d,\Re\langle\cdot,\cdot\rangle)$ such that
    \begin{equation} \label{eq-characteristic-function-of-a-gaussian-state}
        \Tr (\rho W(z)) = \exp{ 2 \mi \Im \langle \omega, z \rangle - \frac{1}{2} \Re \langle z, S z \rangle }, \quad \forall z \in \C^d.
    \end{equation}
        
    The vector $\omega$ and the operator $S$ are called the mean vector and covariance operator, respectively; see \cite{parthasarathy2010gaussian,holevo2011probabilistic} for this characterization. The necessary and sufficient covariance condition is
    \begin{equation*}
        \mathbf{S} + \mi \mathbf{J} \ge 0,
    \end{equation*}
    where $\mathbf S$ is the real matrix representing $S$ (Appendix \ref{section-appendix-matrix-identification}) and
    \begin{equation} \label{eq-boldface-J}
        \mathbf{J} := \begin{bmatrix}
            0 & \1_d \\ 
            - \1_d & 0
        \end{bmatrix}
    \end{equation}
    represents the real-linear operator
    \begin{equation} \label{eq-definition-real-linear-J}
        J : z \mapsto -\mi z.
    \end{equation}
    Thus $J^{-1}=J^\sharp=-J$, where $\cdot^\sharp$ denotes the adjoint with respect to the real inner product $\Re\langle\cdot,\cdot\rangle$.
    
    By Williamson's theorem \cite{williamson1936algebraic}, there is a real symplectic matrix $\mathbf G$ such that
    \begin{equation} \label{eq-symplectic-diagonalization-covariance-matrix}
        \mathbf{S} = \mathbf{G}^* \begin{bmatrix}
            D_{\coth} & 0 \\ 
            0 & D_{\coth}
        \end{bmatrix} \mathbf{G},
    \end{equation}
    where
    \begin{equation} \label{eq-definition-diagonal-matrix-coth}
        D_{\coth} := \diag \left\{ \coth \frac{\beta_1}{2}, \coth \frac{\beta_2}{2}, \cdots, \coth \frac{\beta_d}{2} \right\},
    \end{equation}
    for uniquely determined parameters
    \begin{equation} \label{eq-definition-inverse-temperatures}
        \beta_1,\ldots,\beta_d>0,
    \end{equation}
    up to permutation. We call these parameters the \emph{inverse temperatures}.

    \begin{remark} \label{remark-definition-diagonal-hyperbolic-functions}
        In addition to $D_{\coth}$ in \eqref{eq-definition-diagonal-matrix-coth}, we use
        \[
            D_{\csch}:=\diag\left\{\csch\frac{\beta_1}{2},\ldots,
            \csch\frac{\beta_d}{2}\right\}.
        \]
        Both $D_{\coth}$ and $D_{\csch}$ are also regarded as real-linear operators on $\C^d$, and $\mathbf D_{\coth}$ and $\mathbf D_{\csch}$ denote their $2d\times2d$ real matrix representations. Thus, \eqref{eq-symplectic-diagonalization-covariance-matrix} reads $\mathbf S=\mathbf G^*\mathbf D_{\coth}\mathbf G$.
    \end{remark}
    
    By \cite[Section 4]{parthasarathy2021pedagogical}, a faithful Gaussian state with mean $\omega$ and covariance $S$ can be written as
    \begin{equation} \label{eq-gaussian-state-parthasarathy-form}
        \rho = W(\omega) B(\mathbf{M}) \bigotimes_{j=1}^d (1 - \me^{-\beta_j}) \me^{-\beta_j N_j} B(\mathbf{M})^{-1} W(\omega)^{-1},
    \end{equation}
    where $\mathbf M:=\mathbf G^{-1}$ and $B(\mathbf M)$ is a unitary operator on $\mathsf h$ satisfying
    \begin{equation*}
        B (\mathbf{M}) W (z) B (\mathbf{M})^* = W(M z), \quad \forall z \in \C^d.
    \end{equation*}
    Here $M$ is the real-linear symplectic transformation represented by $\mathbf M$. The unitary $B(\mathbf M)$ exists by the Stone--von Neumann theorem and is unique up to a scalar of modulus one. Such a unitary is called a \emph{Bogoliubov implementer}; see, for example, \cite{bruneau2007bogoliubov}. The real-linear conventions and the decomposition of $M$ into complex-linear and conjugate-linear parts are given in Appendix \ref{section-appendix-matrix-identification}.

    The predual semigroup $(\cT_{*t})_{t\geq0}$ on $\cB_1(\mathsf h)$ is characterized by
    \begin{equation*}
        \Tr \left( \cT_{*t} (\psi) X \right) = \Tr \left( \psi \cT_t (X) \right), \quad \forall \psi \in \cB_1 (\mathsf{h}), \quad \forall X \in \cB (\mathsf{h}).
    \end{equation*}

    The QMS is Gaussian precisely when $(\cT_{*t})_{t\geq0}$ preserves Gaussian states. For the generator in \eqref{eq-gksl-generator}, the predual generator is
    \begin{equation*}
            \cL_*(\psi)=-\mi[H,\psi]-\frac{1}{2}\sum_{\ell=1}^m
            \left(L_\ell^*L_\ell\psi-2L_\ell\psi L_\ell^*+\psi L_\ell^*L_\ell\right),
            \qquad \psi\in\dom\cL_*,
    \end{equation*}
    with the same $H$ and $L_\ell$ as above. The action on Weyl operators completely determines the semigroup and is given by the following formula \cite[Theorem 2.4]{agredo2022decoherence}.

    \begin{theorem} \label{theo-explicit-action-on-weyl-operators}
        For every $t\geq0$ and $z\in\C^d$,
        \begin{equation*}
            \cT_t (W(z)) = \exp{ - \frac{1}{2} \int_0^t \Re \left\langle \me^{s Z} z, C \me^{s Z} z \right\rangle \md s + \mi \int_0^t \Re \left\langle \zeta, \me^{sZ} z \right\rangle \md s} W(\me^{t Z} z),
        \end{equation*}
        where the real-linear operators $Z,C:\C^d\to\C^d$ are
        \begin{align} 
            Z z &= \left( ( U^T \overline{U} - V^T \overline{V} ) / 2 + \mi \Omega \right) z + \left( (U^T V - V^T U)/2 + \mi \kappa \right) \overline{z}, \label{eq-Z-real-operator} \\ 
            C z &= ( U^T \overline{U} + V^T \overline{V} ) z + (U^T V + V^T U) \overline{z}. \label{eq-C-real-operator}
        \end{align}
    \end{theorem}

    We call $Z$ and $C$ the drift and diffusion operators, respectively, and denote their real matrix representations by $\mathbf Z$ and $\mathbf C$.

    If a Gaussian state $\sigma$ has mean $\omega_0$ and covariance $S_0$, then $\cT_{*t}(\sigma)$ has mean $\omega_t$ and covariance $S_t$ given by
    \begin{align}
        \omega_t &= J \me^{t Z^{\sharp}} J^{\sharp} \omega_0 + \frac{1}{2} \int_0^t J \me^{s Z^{\sharp}} \zeta \, \md s, \label{eq-evolution-mean-vector} \\
        S_t &= \me^{t Z^\sharp} S_0 \me^{t Z} + \int_0^t \me^{s Z^{\sharp}} C \me^{s Z} \, \md s, \label{eq-evolution-covariance-operator}
    \end{align}
    where $ J $ is the real-linear operator defined in \eqref{eq-definition-real-linear-J}.

    We call $Z$ (equivalently, $\mathbf Z$) \emph{stable} when every eigenvalue of $\mathbf Z$ has strictly negative real part. Equations \eqref{eq-evolution-mean-vector}--\eqref{eq-evolution-covariance-operator} give the following standard consequence; see also \cite[Theorem 7]{fagnola2025spectral}, after translating conventions.

    \begin{proposition} \label{proposition-Lyapunov-equation}
        If $Z$ is stable, then there is a unique invariant state within the Gaussian class. Its mean $\omega$ and covariance $S$ are
        \begin{equation} \label{eq-relationship-S-C-Z-omega-zeta}
            \omega = \frac{\mi}{2} (Z^\sharp)^{-1} \zeta, \quad S = \int_0^\infty \me^{s Z^\sharp} C \me^{s Z} \, \md s.
        \end{equation}
        In particular, $ S $ satisfies the Lyapunov equation
        \begin{equation} \label{eq-Lyapunov-equation}
            Z^\sharp S + S Z = - C.
        \end{equation} 
    \end{proposition}

    Proposition \ref{proposition-Lyapunov-equation} establishes uniqueness only within the Gaussian class. We impose the following standing assumption.

    \begin{assumption*}[Standing assumption]
        The Gaussian QMS $(\cT_t)_{t\geq0}$ has a unique faithful normal invariant state $\rho$.
    \end{assumption*}

    This assumption is equivalent to the conjunction of stability and strict covariance positivity stated below; see \cite[proof of Proposition 15]{GIROTTI2026130150}. It is also related to irreducibility; see \cite{fagnola2026irreducibility}.

    \begin{proposition} \label{proposition-unique-invariant-faithful-state}
        The Gaussian QMS $(\cT_t)_{t\geq0}$ has a unique faithful normal invariant state if and only if $\mathbf Z$ is stable and $\mathbf S+\mi\mathbf J>0$, where $S$ is the covariance operator in \eqref{eq-relationship-S-C-Z-omega-zeta}.
    \end{proposition}
    
    Consequently, $\rho$ has the form \eqref{eq-gaussian-state-parthasarathy-form}. Using Williamson's theorem \cite{williamson1936algebraic} and the unitary implementation in \cite[Section 4]{parthasarathy2021pedagogical}, we conjugate the semigroup by the unitary appearing there. This gives a unitarily equivalent Gaussian QMS whose invariant state is
    \begin{equation} \label{eq-diagonalized-gaussian-state}
        \rho_0 := \bigotimes_{j=1}^d (1 - \me^{-\beta_j}) \me^{-\beta_j N_j}.
    \end{equation}
    We call this procedure \emph{standardization}. Appendix \ref{section-appendix-invariance-spectral-gap} proves that it yields a unitary equivalence of the induced semigroups and therefore preserves their spectral gaps.
    
    The state $\rho_0$ has mean zero and covariance $D_{\coth}$. Until the end of Section \ref{section-application-boson-chain-model}, we work with the standardized semigroup and suppress tildes from its parameters.

    \section{Existence of the KMS Spectral Gap} \label{section-existence-kms-spectral-gap}
    
    By \cite[Proposition 14 and Theorem 27]{fagnola2025spectral}, the GNS spectral gap is $-1/2$ times the largest eigenvalue of
    \begin{equation} \label{eq-gns-spectral-gap-matrix}
        \mathbf{Z} + ( \mathbf{D}_{\coth} + \mi \mathbf{J} )^{-1} \mathbf{Z}^* ( \mathbf{D}_{\coth} + \mi \mathbf{J} ),
    \end{equation} 
    whereas the KMS spectral gap is $-1/2$ times the largest eigenvalue of
    \begin{equation} \label{eq-kms-spectral-gap-matrix} 
        \mathbf{Z} + \mathbf{D}_{\csch}^{-1} \mathbf{Z}^* \mathbf{D}_{\csch}.
    \end{equation}
    The matrices in \eqref{eq-gns-spectral-gap-matrix} and \eqref{eq-kms-spectral-gap-matrix} are similar to self-adjoint, negative semidefinite matrices \cite[Proposition 14 and Section 7]{fagnola2025spectral}; in particular, their spectra are real and nonpositive. For the KMS matrix, the similarity is explicit:
    \[
        \mathbf D_{\csch}^{1/2}
        \bigl(\mathbf Z+\mathbf D_{\csch}^{-1}\mathbf Z^*\mathbf D_{\csch}\bigr)
        \mathbf D_{\csch}^{-1/2}
        =\mathbf D_{\csch}^{-1/2}
        \bigl(\mathbf Z^*\mathbf D_{\csch}+\mathbf D_{\csch}\mathbf Z\bigr)
        \mathbf D_{\csch}^{-1/2},
    \]
    and the matrix on the right is self-adjoint and negative semidefinite. Moreover, the GNS matrix in \eqref{eq-gns-spectral-gap-matrix} has strictly negative spectrum precisely when $[U\,\overline V]$ has full column rank \cite[Proposition 6, Theorem 18, and Proposition 22]{fagnola2025spectral}. We now derive the corresponding, weaker criterion for the KMS embedding.

    Since $\mathbf D_{\csch}$ is positive definite, positivity of the KMS spectral gap is equivalent to invertibility of the negative semidefinite matrix
    \begin{equation} \label{eq-kms-symmetric-matrix}
        \mathbf{Z}^* \mathbf{D}_{\csch} + \mathbf{D}_{\csch} \mathbf{Z}.
    \end{equation}
    The next identity expresses the corresponding real-linear operator $Z^\sharp D_{\csch}+D_{\csch}Z$ entirely in terms of the Kraus coefficients and makes its negativity transparent.

    \begin{theorem}[KMS dissipation identity] \label{main-theorem-1}
        For every $z\in\C^d$,
        \begin{align}
            (Z^\sharp D_{\csch} + D_{\csch} Z) z &= \left[  - \frac{2 (U^T \overline{U})_{j k}}{\me^{-\beta_j/2} + \me^{-\beta_k/2}}  - \frac{2 (V^T \overline{V})_{j k}}{\me^{\beta_j/2} + \me^{\beta_k/2}}  \right]_{j, k} z \nonumber \\ 
            &\quad + \left[ - \frac{2 (U^T V)_{j k} }{\me^{-\beta_j/2} + \me^{\beta_k/2}}  - \frac{2 (V^T U)_{j k}}{\me^{\beta_j/2} + \me^{-\beta_k/2}} \right]_{j, k} \overline{z}, \label{eq-explicit-Z-sharp-Dcsch}
        \end{align}
        where each bracket denotes the $d\times d$ matrix with indices $1\leq j,k\leq d$. In particular,
        \begin{equation} \label{eq-KMS-dissipation-integral}
            Z^\sharp D_{\csch}+D_{\csch}Z
            =-2\int_0^\infty A_t^\sharp A_t\,\md t\leq0,
        \end{equation}
        where $A_t:\C^d\to\C^m$ is the real-linear map
        \begin{equation} \label{eq-definition-A_t-real-linear-op}
            A_t z := \left[ \exp{- \me^{-\beta_k / 2} t} \overline{u_{\ell k}} \right]_{\ell, k} z + \left[ \exp{-\me^{\beta_k / 2} t} v_{\ell k} \right]_{\ell, k} \overline{z},
        \end{equation}
        with $1\leq \ell\leq m$ and $1\leq k\leq d$.
    \end{theorem}
    \begin{proof}
        The Lyapunov equation \eqref{eq-Lyapunov-equation} becomes
        \begin{equation} \label{eq-Lyapunov-equation-diagonal-by-coth}
            Z^\sharp D_{\coth} + D_{\coth} Z = - C.
        \end{equation}
        Equating its complex-linear and conjugate-linear parts and using \eqref{eq-Z-real-operator}--\eqref{eq-C-real-operator}, we obtain
        \begin{align}
            \left\{ (U^T \overline{U} - V^T \overline{V})/2, D_{\coth} \right\} + \mi \left[ D_{\coth}, \Omega \right] = -( U^T \overline{U} + V^T \overline{V}), \label{eq-lyapunov-equation-by-parameter-complex-linear-part} \\
            \left[ (V^T U - U^T V)/2, D_{\coth} \right] + \mi \left\{ \kappa, D_{\coth} \right\} = -(U^T V + V^T U). \label{eq-lyapunov-equation-by-parameter-anti-linear-part}
        \end{align}
        Taking the $(j,k)$-entries gives
        \begin{align}
            \mi \Omega_{j k} \coth(\beta_j /2) - \mi \Omega_{j k} \coth(\beta_k / 2) &= - ( U^T \overline{U} + V^T \overline{V} )_{j k} - \frac{1}{2} ( U^T \overline{U} - V^T \overline{V} )_{j k} \coth(\beta_k /2) \nonumber \\ 
            &\qquad - \frac{1}{2} (U^T \overline{U} - V^T \overline{V})_{j k} \coth(\beta_j / 2) \label{eq-lyapunov-equation-by-parameter-complex-linear-part-j-k-term} \\ 
            \mi \kappa_{j k} \coth(\beta_j/2) + \mi \kappa_{j k} \coth(\beta_k/2) &= -( U^T V + V^T U )_{j k} - \frac{1}{2} (V^T U - U^T V)_{j k} \coth(\beta_k/2) \nonumber \\ 
            &\qquad -\frac{1}{2} (U^T V - V^T U)_{jk} \coth(\beta_j / 2), \label{eq-lyapunov-equation-by-parameter-anti-linear-part-j-k-term}
        \end{align}
        which implies that, for $ \beta_j \neq \beta_k $,
        \begin{align}
            \mi \Omega_{j k} &= \frac{ -(\coth(\beta_j/2) + \coth(\beta_k/2) )/2 - 1}{\coth(\beta_j/2) - \coth(\beta_k/2)} (U^T \overline{U})_{j k} \nonumber \\
            &\quad + \frac{ (\coth(\beta_j/2) + \coth(\beta_k/2) )/2 - 1}{\coth(\beta_j/2) - \coth(\beta_k/2)} (V^T \overline{V})_{j k} \nonumber \\
            &= \frac{2 \me^{\beta_k} - 1 - \me^{\beta_k} \me^{-\beta_j}}{2 (1 - \me^{-\beta_j} \me^{\beta_k})} (U^T \overline{U})_{j k} + \frac{2 \me^{-\beta_j} - 1 - \me^{\beta_k} \me^{-\beta_j}}{2 (1 - \me^{-\beta_j} \me^{\beta_k})} (V^T \overline{V})_{j k}, \label{eq-entries-i-Omega-jk} 
        \end{align}
        and, for all $1\leq j,k\leq d$,
        \begin{align}
            \mi \kappa_{j k} &= \frac{ (\coth(\beta_k/2) - \coth(\beta_j/2) )/2 - 1}{\coth(\beta_k/2) + \coth(\beta_j/2)} (U^T V)_{j k} \nonumber \\ 
            &\quad + \frac{ (\coth(\beta_j/2) - \coth(\beta_k/2) )/2 - 1}{\coth(\beta_k/2) + \coth(\beta_j/2)} (V^T U)_{j k} \nonumber \\ 
            &= \frac{2 \me^{\beta_j} - \me^{\beta_k} \me^{\beta_j} - 1}{2 (\me^{\beta_k} \me^{\beta_j} - 1)} (U^T V)_{j k} + \frac{2 \me^{\beta_k} - \me^{\beta_k} \me^{\beta_j} - 1}{2 (\me^{\beta_k} \me^{\beta_j} - 1)} (V^T U)_{j k}. \label{eq-entries-i-kappa-jk}
        \end{align}
       
        When $\beta_j=\beta_k$, the complex-linear part of the Lyapunov equation does not determine $\Omega_{jk}$. This entry is irrelevant to $Z^\sharp D_{\csch}+D_{\csch}Z$, because the corresponding commutator with $D_{\csch}$ vanishes. In this case, \eqref{eq-lyapunov-equation-by-parameter-complex-linear-part-j-k-term} instead gives
        \begin{equation*}
            - (U^T \overline{U} + V^T \overline{V})_{j k} - (U^T \overline{U} - V^T \overline{V})_{j k} \coth(\beta_j / 2) = 0,
        \end{equation*}
        and then
        \begin{equation} \label{eq-equal-temperature-UU-VV}
            (U^T \overline{U})_{j k } = \me^{-\beta_j} ( V^T \overline{V} )_{j k}.
        \end{equation}

        From the definition \eqref{eq-Z-real-operator},
        \begin{align}
            (Z^\sharp D_{\csch} + D_{\csch} Z)z &= \left( ( U^T \overline{U} - V^T \overline{V} ) / 2 - \mi \Omega \right) D_{\csch} z + \left( (V^T U - U^T V ) / 2 + \mi \kappa \right) D_{\csch} \overline{z}  \nonumber \\ 
            &\quad + D_{\csch} \left( (U^T \overline{U} - V^T \overline{V}) / 2  + \mi \Omega \right) z + D_{\csch} \left( (U^T V - V^T U) / 2 + \mi \kappa \right) \overline{z}  \nonumber \\ 
            &= \left(  \left\{ ( U^T \overline{U} - V^T \overline{V} ) / 2, D_{\csch} \right\}  +  \mi [D_{\csch}, \Omega] \right) z  \nonumber \\  
            &\quad + \left( \left[ (V^T U - U^T V) / 2, D_{\csch} \right] + \mi \left\{ \kappa, D_{\csch} \right\} \right) \overline{z}, \label{eq-explicit-form-Z-sharp-D-D-Z}
        \end{align}

        Substituting \eqref{eq-entries-i-Omega-jk} and \eqref{eq-entries-i-kappa-jk} into the $(j,k)$-entries of \eqref{eq-explicit-form-Z-sharp-D-D-Z}, we obtain the following. If $\beta_j\neq\beta_k$, then
        \begin{align}
            &\quad \left( \left\{ (U^T \overline{U} - V^T \overline{V})/2, D_{\csch} \right\} + \mi \left[ D_{\csch}, \Omega \right] \right)_{jk} \nonumber \\
            &= \frac{\csch(\beta_k /2) + \csch(\beta_j / 2)}{2} \left( (U^T \overline{U})_{j k} - (V^T \overline{V})_{j k} \right) \nonumber \\ 
            &\quad + \mi \Omega_{j k} \left( \csch(\beta_j / 2) - \csch (\beta_k / 2) \right) \nonumber \\
            &= \left( \frac{\me^{\beta_k/2}}{\me^{\beta_k} - 1} + \frac{\me^{\beta_{j}/2}}{\me^{\beta_j} - 1} + \left( \frac{\me^{\beta_j / 2}}{\me^{\beta_j} - 1} - \frac{\me^{\beta_k/2}}{\me^{\beta_k} - 1} \right) \frac{2 \me^{\beta_j} \me^{\beta_k} - \me^{\beta_j} - \me^{\beta_k}}{\me^{\beta_j} - \me^{\beta_k}} \right) (U^T \overline{U})_{jk} \nonumber \\
            &\quad + \left( - \frac{\me^{\beta_k / 2}}{\me^{\beta_k} - 1} - \frac{\me^{\beta_j / 2}}{\me^{\beta_j} - 1} + \left( \frac{\me^{\beta_j / 2}}{\me^{\beta_j} - 1} - \frac{\me^{\beta_k / 2}}{\me^{\beta_k} - 1} \right) \frac{2 - \me^{\beta_j} - \me^{-\beta_k}}{\me^{\beta_j} - \me^{\beta_k}} \right) (V^T \overline{V})_{j k} \nonumber \\
            &= - \frac{2 (U^T \overline{U})_{j k}}{\me^{-\beta_j/2} + \me^{-\beta_k/2}}  - \frac{2 (V^T \overline{V})_{j k}}{\me^{\beta_j/2} + \me^{\beta_k/2}}, \label{eq-explicit-form-Z-sharp-D-D-Z-complex-linear-part}
        \end{align}
        If $\beta_j=\beta_k$, then \eqref{eq-equal-temperature-UU-VV} gives
        \begin{align*}
            &\quad \left( \left\{ (U^T \overline{U} - V^T \overline{V})/2, D_{\csch} \right\} + \mi \left[ D_{\csch}, \Omega \right] \right)_{jk} \\
            &= \frac{\csch(\beta_k /2) + \csch(\beta_j / 2)}{2} \left( (U^T \overline{U})_{j k} - (V^T \overline{V})_{j k} \right) \\ 
            &= - 2 \me^{-\beta_j/2} ( V^T \overline{V} )_{j k} = - \frac{2 (U^T \overline{U})_{j k}}{\me^{-\beta_j/2} + \me^{-\beta_k/2}}  - \frac{2 (V^T \overline{V})_{j k}}{\me^{\beta_j/2} + \me^{\beta_k/2}}.
        \end{align*}
        Thus \eqref{eq-explicit-form-Z-sharp-D-D-Z-complex-linear-part} also holds when $\beta_j=\beta_k$.
        
        Moreover,
        \begin{align}
            &\quad \left( \left[ (V^T U - U^T V)/2, D_{\csch} \right] + \mi \left\{ \kappa, D_{\csch} \right\} \right)_{jk} \nonumber \\
            &= \frac{\csch(\beta_j /2) - \csch(\beta_k / 2)}{2} \left( (U^T V)_{j k} - (V^T U)_{j k} \right) \nonumber \\ 
            &\quad + \mi \kappa_{j k} \left( \csch(\beta_j / 2) + \csch (\beta_k / 2) \right) \nonumber \\
            &= \left( \frac{\me^{\beta_j / 2}}{\me^{\beta_j} - 1} - \frac{\me^{\beta_k / 2}}{\me^{\beta_k} - 1} + \left( \frac{\me^{\beta_j / 2}}{\me^{\beta_j} - 1} + \frac{\me^{\beta_k / 2}}{\me^{\beta_k} - 1}  \right) \frac{2 \me^{\beta_j} - \me^{\beta_j} \me^{\beta_k} - 1}{\me^{\beta_j} \me^{\beta_k} - 1} \right) (U^T V)_{j k} \nonumber \\
            &\quad + \left( \frac{\me^{\beta_k / 2}}{\me^{\beta_k} - 1} - \frac{\me^{\beta_j / 2}}{\me^{\beta_j} - 1} + \left( \frac{\me^{\beta_j / 2}}{\me^{\beta_j} - 1} + \frac{\me^{\beta_k / 2}}{\me^{\beta_k} - 1} \right) \frac{2 \me^{\beta_k} - \me^{\beta_k} \me^{\beta_j} - 1}{\me^{\beta_k} \me^{\beta_j} - 1} \right) (V^T U)_{jk} \nonumber \\
            &= - \frac{2 (U^T V)_{j k} }{\me^{-\beta_j/2} + \me^{\beta_k/2}}  - \frac{2 (V^T U)_{j k}}{\me^{\beta_j/2} + \me^{-\beta_k/2}}. \label{eq-explicit-form-Z-sharp-D-D-Z-anti-linear-part}
        \end{align}
        These two identities give the complex-linear and conjugate-linear parts of $Z^\sharp D_{\csch}+D_{\csch}Z$, respectively. This proves \eqref{eq-explicit-Z-sharp-Dcsch}.

        Now, by \eqref{eq-definition-A_t-real-linear-op}, for $y\in\C^m$ we have
        \begin{equation*}
            A_t^\sharp y = \left[ \exp{-\me^{-\beta_j / 2} t} u_{\ell j} \right]_{j, \ell} y + \left[ \exp{-\me^{\beta_j / 2} t} v_{\ell j} \right]_{j, \ell} \overline{y},
        \end{equation*}
        where $1\leq j\leq d$ and $1\leq \ell\leq m$. Integrating the products of the exponential factors gives \eqref{eq-KMS-dissipation-integral}.
    \end{proof}

    We now group the columns of $U$ and $V$ according to the inverse temperatures in \eqref{eq-definition-diagonal-matrix-coth}. Some of the inverse temperatures $\beta_1,\ldots,\beta_d$ may coincide. Let $\beta_{(1)},\ldots,\beta_{(l)}$ denote their distinct values. For $1\leq n\leq l$, define
    \begin{equation} \label{eq-definition-I-n}
        I_n := \left\{ j \in \{1,\ldots,d\} : \beta_j = \beta_{(n)} \right\}.
    \end{equation}
    Thus, $I_1,\ldots,I_l$ form a partition of $\{1,\ldots,d\}$. We call each $I_n$ an \emph{equal-temperature block}. With the indices in $I_n$ arranged increasingly, define the corresponding column submatrices
    \begin{equation} \label{eq-definition-U-V-n}
        U_n := \begin{bmatrix}
            U e_i
        \end{bmatrix}_{i \in I_n}, \quad V_n := \begin{bmatrix}
            V e_i
        \end{bmatrix}_{i \in I_n},
    \end{equation}
    where $e_i$ is the $i$-th canonical vector of $\C^d$. Thus $U_n,V_n\in M_{m\times |I_n|}(\C)$, and their kernels are subspaces of $\C^{|I_n|}$. Our notation is fixed throughout: $d$ is the number of modes, $m$ is the number of Kraus operators, $l$ is the number of equal-temperature blocks, and $\ell$ is the index of a Kraus operator.

    The equal-temperature part of the Lyapunov equation gives the following basic relation.

    \begin{lemma} \label{lemma-same-block-relation-U-and-V}
        For every $1\leq n\leq l$,
        \begin{equation} \label{eq-block-relation-U-and-V}
            U^*_n U_n = \me^{-\beta_{(n)}} V_n^* V_n, 
        \end{equation}
        Consequently,
        \begin{equation} \label{eq-identical-kernels-in-block}
            \ker U_n = \ker V_n.
        \end{equation}
    \end{lemma}
    \begin{proof}
        Applying \eqref{eq-equal-temperature-UU-VV} to all pairs of indices in $I_n$ gives
        \begin{equation*}
            U^T_n \overline{U}_n = \me^{-\beta_{(n)}} V^T_n \overline{V}_n.
        \end{equation*}
        Taking entrywise complex conjugates yields \eqref{eq-block-relation-U-and-V}. Since $\ker A=\ker(A^*A)$ for every complex matrix $A$,
        \begin{equation*}
            \ker U_n = \ker ( U_n^* U_n ) = \ker ( V_n^* V_n ) = \ker V_n,
        \end{equation*} 
        which gives \eqref{eq-identical-kernels-in-block}.
    \end{proof}
    
    We can now state the main theorem.

    \begin{theorem}[Main theorem: blockwise rank criterion] \label{theorem-main-2}
        The following are equivalent:
        \begin{enumerate}
            \item The KMS spectral gap $\gap L^{(1/2)}$ is strictly positive. \label{item-1}
            \item $\ker U_n=\{0\}$ for every $1\leq n\leq l$. \label{item-2}
            \item $\ker V_n=\{0\}$ for every $1\leq n\leq l$. \label{item-3}
            \item $\rank U_n=\rank V_n=|I_n|$ for every $1\leq n\leq l$. \label{item-4}
        \end{enumerate}
    \end{theorem}
    \begin{proof}
        By Lemma \ref{lemma-same-block-relation-U-and-V}, conditions \eqref{item-2}--\eqref{item-4} are equivalent. It remains to characterize condition \eqref{item-1}.

        Set $K:=Z^\sharp D_{\csch}+D_{\csch}Z$. The discussion preceding Theorem \ref{main-theorem-1} and the negativity of $K$ show that $\gap L^{(1/2)}>0$ if and only if $\ker K=\{0\}$.
        By \eqref{eq-KMS-dissipation-integral}, for every $w\in\C^d$,
        \begin{equation} \label{eq-main-theorem-quadratic-form}
            -\Re\langle w,Kw\rangle
            =2\int_0^\infty \norm{A_t w}^2\,\md t.
        \end{equation}
        Since $t\mapsto A_t w$ is continuous, $w\in\ker K$ if and only if $A_t w=0$ for every $t\geq0$.

        Write $w_{(n)}=(w_k)_{k\in I_n}\in\C^{|I_n|}$. By \eqref{eq-definition-A_t-real-linear-op}, the identity $A_t w=0$ for all $t\geq0$ is equivalent, componentwise, to
        \[
            \sum_{n=1}^l \left(
            \exp{-\me^{-\beta_{(n)}/2}t}\,\overline{U_n}w_{(n)}
            +\exp{-\me^{\beta_{(n)}/2}t}\,V_n\overline{w_{(n)}}
            \right)=0, \qquad t\geq0.
        \]
        The $2l$ positive rates $\me^{-\beta_{(n)}/2}$ and $\me^{\beta_{(n)}/2}$ are pairwise distinct. Within each family they are distinct because different blocks have different inverse temperatures; moreover, the first family lies in $(0,1)$, whereas the second lies in $(1,\infty)$, so no rate from one family can equal a rate from the other. The functions $t\mapsto\exp{-\lambda t}$ with distinct rates $\lambda$ are linearly independent. Hence, for every $n$,
        \[
            \overline{U_n}w_{(n)}=0,
            \qquad V_n\overline{w_{(n)}}=0.
        \]
        Taking the complex conjugate of the first relation gives $U_n\overline{w_{(n)}}=0$, while the second relation already gives $V_n\overline{w_{(n)}}=0$. Thus $\overline{w_{(n)}}\in\ker U_n\cap\ker V_n$. Hence $\ker K=\{0\}$ if and only if $\ker U_n\cap\ker V_n=\{0\}$ for every $n$.
        Lemma \ref{lemma-same-block-relation-U-and-V} identifies the two kernels in each block, completing the proof.
    \end{proof}

    \begin{corollary} \label{corollary-GNS-implies-KMS}
        If the GNS spectral gap is positive, then the KMS spectral gap is positive.
    \end{corollary}
    \begin{proof}
        Positivity of the GNS spectral gap is equivalent to
        $\rank[U\,\overline V]=2d$. Hence $U$ and $\overline V$, and therefore $U$ and $V$, have full column rank. Every column submatrix $U_n$ and $V_n$ is therefore injective, so
        \begin{equation*}
            \ker U_n=\ker V_n=\{0\}, \qquad 1\leq n\leq l.
        \end{equation*}
        The conclusion follows from Theorem \ref{theorem-main-2}.
    \end{proof}

    \begin{remark}[Comparison with the GNS criterion]
        The GNS criterion requires global linear independence of all $2d$ columns of $[U\,\overline V]$. The KMS criterion requires only that the columns of $U$ (equivalently, of $V$) be linearly independent within each equal-temperature block. It does not require independence between columns belonging to different blocks, nor between the $U$- and $V$-columns.
    \end{remark}

    \begin{remark}[Comparison of the gap constants]
        It was conjectured in \cite[Section 7]{fagnola2025spectral} that the KMS spectral gap dominates the GNS spectral gap. Wirth proved this inequality for every QMS with a faithful normal invariant state \cite[Corollary 3.7]{wirth2026kms}. Corollary \ref{corollary-GNS-implies-KMS} is therefore also a consequence of that general result. The proof above is retained because it directly compares the two Gaussian rank conditions.
    \end{remark}

    The block structure makes the role of the inverse temperatures explicit. If all inverse temperatures are distinct, every block is a singleton and the KMS criterion reduces to the nonvanishing of each corresponding column. The next lemma shows that the standing assumption rules out a block on which both Kraus-coefficient matrices vanish.

    \begin{lemma} \label{lemma-zero-Un-Vn-not-possible}
        For every $1\leq n\leq l$, one has $U_n\neq0$ and $V_n\neq0$.
    \end{lemma}
    \begin{proof}
        Suppose that $U_n=0$. Lemma \ref{lemma-same-block-relation-U-and-V} then gives $V_n=0$. Thus, for each $j\in I_n$, every entry involving the $j$-th column vanishes in $U^T\overline U$, $V^T\overline V$, $U^TV$, and $V^TU$. For $1\leq k\leq d$, the $(j,k)$-entries of \eqref{eq-lyapunov-equation-by-parameter-complex-linear-part-j-k-term} and \eqref{eq-lyapunov-equation-by-parameter-anti-linear-part-j-k-term} reduce to
        \begin{equation*}
            \mi \Omega_{j k} ( \coth(\beta_j / 2) - \coth(\beta_k /2 ) ) = 0, \quad \mi \kappa_{j k} ( \coth(\beta_j / 2) + \coth(\beta_k /2 ) ) = 0.
        \end{equation*}
        Since every $\beta_k$ is positive, it follows that
        \begin{equation} \label{eq-entries-Omega-kappa-when-zero-column}
            \Omega_{j k} = 0, \quad \forall k \notin I_n; \quad \kappa_{j k} = 0, \quad \forall 1 \le k \le d.
        \end{equation}
        Define $E_n:=\myspan\{e_j:j\in I_n\}$. Since $\Omega=\Omega^*$ and $\kappa=\kappa^T$, \eqref{eq-entries-Omega-kappa-when-zero-column} implies
        \begin{equation*}
            \Omega E_n \subset E_n, \quad \kappa E_n = \{ 0 \}.
        \end{equation*}
        Moreover, $U\restriction_{E_n}=V\restriction_{E_n}=0$. Hence \eqref{eq-Z-real-operator} gives
        \begin{equation*}
            Z\restriction_{E_n}=\mi\,\Omega\restriction_{E_n}.
        \end{equation*}
        The restriction $\Omega\restriction_{E_n}$ is Hermitian, so $Z\restriction_{E_n}=\mi\,\Omega\restriction_{E_n}$ has spectrum on the imaginary axis. This contradicts the stability of $Z$. Thus $U_n\neq0$, and Lemma \ref{lemma-same-block-relation-U-and-V} also gives $V_n\neq0$.
    \end{proof}

    \begin{corollary} \label{corollary-distinct-temperatures}
        If the inverse temperatures $\beta_1,\ldots,\beta_d$ are pairwise distinct, then the KMS spectral gap is positive. In particular, this holds when $d=1$.
    \end{corollary}
    \begin{proof}
        Every equal-temperature block is a singleton. Lemma \ref{lemma-zero-Un-Vn-not-possible} shows that its unique $U$-column is nonzero, and Theorem \ref{theorem-main-2} applies.
    \end{proof}

    Repeated inverse temperatures allow a different phenomenon: a block may be nonzero but fail to have full column rank. The following two-mode example satisfies the standing assumption but has no KMS spectral gap.

    \begin{example} \label{example-two-modes-the-same-inverse-temperatures}
        Let $d=2$, choose $v>u>0$, and set $U=[u\ u]$ and $V=[v\ v]$. There is one Kraus operator,
        \[
            L=v a_1+u a_1^\dagger+v a_2+u a_2^\dagger,
        \]
        and
        \begin{equation*}
            U^T V = V^T U = v u \begin{bmatrix}
                1 & 1 \\ 
                1 & 1
            \end{bmatrix}, \quad U^T \overline{U} = u^2 \begin{bmatrix}
                1 & 1 \\ 
                1 & 1
            \end{bmatrix}, \quad V^T \overline{V} = v^2 \begin{bmatrix}
                1 & 1 \\ 
                1 & 1
            \end{bmatrix}.
        \end{equation*}
        Choose
        \begin{equation*}
            \kappa = \mi u v \cdot \frac{v^2 - u^2}{v^2 + u^2} \begin{bmatrix}
                1 & 1 \\ 
                1 & 1
            \end{bmatrix}, \quad \Omega = \begin{bmatrix}
            a & b \\ b & c
            \end{bmatrix},
        \end{equation*}
        where $a,b,c\in\R$ and $a\neq c$. Let
        $E:=\begin{bmatrix}1&1\\1&1\end{bmatrix}$. The drift operator is
        \begin{equation*}
            Z z = \left( \frac{u^2 - v^2}{2} E + \mi \Omega \right) z + \left( - u v \cdot \frac{v^2-u^2}{v^2 + u^2} E \right) \overline{z},
        \end{equation*}
        and its real matrix representation is
        \begin{equation} \label{eq-example-Z}
            \mathbf{Z} = \begin{bmatrix}
                \left( \frac{u^2 - v^2}{2} - u v \cdot \frac{v^2 - u^2}{v^2 + u^2}  \right) E & - \Omega \\ 
                \Omega & \left( \frac{u^2 - v^2}{2} + u v \cdot \frac{v^2 - u^2}{v^2 + u^2}  \right) E
            \end{bmatrix}.
        \end{equation}
        Its Hermitian part is
        \begin{equation} \label{eq-example-Z-star-plus-Z}
            \frac{1}{2} ( \mathbf{Z}^* + \mathbf{Z} ) = \begin{bmatrix}
                \left( \frac{u^2 - v^2}{2} - u v \cdot \frac{v^2 - u^2}{v^2 + u^2}  \right) E & 0 \\ 
                0 & \left( \frac{u^2 - v^2}{2} + u v \cdot \frac{v^2 - u^2}{v^2 + u^2}  \right) E
            \end{bmatrix}.
        \end{equation}

        Since $v>u>0$, the matrix in \eqref{eq-example-Z-star-plus-Z} is negative semidefinite and its kernel is
        \[
            \myspan\{(1,-1,0,0),(0,0,1,-1)\}.
        \]
        If $\mathbf Zw=\lambda w$, then
        \[
            \Re\lambda\,\norm{w}^2
            =\left\langle w,\frac{\mathbf Z^*+\mathbf Z}{2}w\right\rangle\leq0.
        \]
        Equality would force $w=(x,-x,y,-y)$ for some $x,y\in\C$. For such a vector, the condition that $\mathbf Zw$ again belong to this two-dimensional kernel gives $(a-c)x=(a-c)y=0$. Since $a\neq c$, no nonzero eigenvector can satisfy $\Re\lambda=0$. Hence $\mathbf Z$ is stable.

        Proposition \ref{proposition-Lyapunov-equation} now gives a unique invariant Gaussian state. A direct computation shows that its covariance operator is
        \begin{equation*}
            S z = \begin{bmatrix}
                \frac{v^2 + u^2}{v^2 - u^2} & 0 \\
                0 & \frac{v^2 + u^2}{v^2 - u^2}
            \end{bmatrix} z
        \end{equation*}
        because this operator solves \eqref{eq-Lyapunov-equation}. Since $(v^2+u^2)/(v^2-u^2)>1$, Proposition \ref{proposition-unique-invariant-faithful-state} shows that the invariant state is faithful and is the unique normal invariant state. The two modes have the same inverse temperature,
        \begin{equation} \label{eq-example-same-inverse-temperatures}
            \beta_1 = \beta_2 = \log \left( v^2 / u^2 \right),
        \end{equation}
        so there is a single equal-temperature block $I_1=\{1,2\}$. Since $U_1=[u\ u]$ has rank one, Theorem \ref{theorem-main-2} shows that the KMS spectral gap vanishes. Equivalently, $D_{\csch}$ is a scalar multiple of the identity and
        \begin{equation*}
            Z^\sharp D_{\csch} + D_{\csch} Z = \frac{2 u v}{v^2 - u^2} ( Z^\sharp + Z ),
        \end{equation*}
        whose kernel is nontrivial by \eqref{eq-example-Z-star-plus-Z}.
    \end{example}

    \section{Application: A Three-Mode Chain} \label{section-application-boson-chain-model}

    We apply Corollary \ref{corollary-distinct-temperatures} to a three-mode example. Fix $\omega\in\R\setminus\{0\}$ and two parameters $\theta_1,\theta_3>0$. Let
    \begin{equation} \label{eq-boson-chain-model-Hamiltonian}
        H = \omega \sum_{j=1}^{2}\left(a_j a^\dagger_{j+1}+ a^\dagger_j a_{j+1} \right),
    \end{equation}
    and
    \begin{equation} \label{eq-boson-chain-model-Kraus-operators}
        L_1 = \gamma_1^{-} a_1, \quad L_2 = \gamma_1^+ a_1^\dagger, \quad L_3 = \gamma_3^- a_3, \quad L_4 = \gamma_3^+ a_3^\dagger,
    \end{equation}
    where
    \begin{equation*}
        \gamma_j^- = \sqrt{\frac{\me^{\theta_j}}{\me^{\theta_j} - 1}}, \quad
        \gamma_j^+ = \sqrt{\frac{1}{\me^{\theta_j} - 1}},
        \qquad j\in\{1,3\}.
    \end{equation*}

    The matrix $[U\,\overline V]$ has at most four independent rows and therefore cannot have rank $2d=6$; hence the GNS spectral gap vanishes. We assume from now on that $\theta_1\neq\theta_3$. The parameters $\theta_1$ and $\theta_3$ occur in the coefficients of the boundary Kraus operators; they are not the inverse temperatures obtained after Williamson diagonalization of the invariant covariance matrix.

    Equations \eqref{eq-boson-chain-model-Hamiltonian}--\eqref{eq-boson-chain-model-Kraus-operators} give
    \begin{equation*}
        U = \begin{bmatrix}
            0 & 0 & 0 \\ 
            \gamma_1^+ & 0 & 0\\ 
            0 & 0 & 0 \\ 
            0 & 0 & \gamma_3^+
        \end{bmatrix}, \quad V = \begin{bmatrix}
            \gamma_1^- & 0 & 0 \\ 
            0 & 0 & 0 \\ 
            0 & 0 & \gamma_3^- \\ 
            0 & 0 & 0
        \end{bmatrix}, \quad \Omega = \begin{bmatrix}
            0 & \omega & 0 \\ 
            \omega & 0 & \omega \\ 
            0 & \omega & 0
        \end{bmatrix}.
    \end{equation*}
    Here $\kappa=0$, $\zeta=0$, $U^TV=V^TU=0$, and
    \begin{align*}
        U^T \overline{U} - V^T \overline{V} &= \diag \{ -1, 0, -1 \}, \\ 
        U^T \overline{U} + V^T \overline{V} &= \diag \{ \coth(\theta_1/2), 0, \coth(\theta_3 / 2) \}.
    \end{align*}

    Therefore,
    \begin{equation*}
        \mathbf{Z} = \begin{bmatrix}
            - \frac{1}{2} \diag \{ 1, 0, 1 \} & - \Omega \\ 
            \Omega &  - \frac{1}{2} \diag \{ 1, 0, 1 \}
        \end{bmatrix}
    \end{equation*}
    and 
    \begin{equation*}
        \mathbf{C} = \begin{bmatrix}
            \diag \{ \coth (\theta_1 / 2), 0, \coth (\theta_3 / 2) \} & 0 \\ 
            0 & \diag \{ \coth (\theta_1 / 2), 0, \coth (\theta_3 / 2) \}
        \end{bmatrix}.
    \end{equation*}

    The Hermitian part of $\mathbf Z$ is negative semidefinite, with kernel corresponding to vectors supported on the middle mode. If $\mathbf Zw=\lambda w$ and $\Re\lambda=0$, then $w$ must belong to this kernel; the coupling $\omega\neq0$ then shows that $\mathbf Zw$ cannot be a scalar multiple of $w$ unless $w=0$. Thus $\mathbf Z$ is stable, by the same numerical-range argument used in Example \ref{example-two-modes-the-same-inverse-temperatures}.

    Following the decomposition used in \cite[Appendix A]{dhahri2024environment}, write
    \begin{equation} \label{eq-boson-chain-model-decomposition-C-matrix}
        \mathbf{C} = \frac{1}{2} \left( \coth \frac{\theta_1}{2} + \coth \frac{\theta_3}{2} \right) \mathbf{C}_1 + \frac{1}{2} \left( \coth \frac{\theta_1}{2} - \coth \frac{\theta_3}{2} \right) \mathbf{C}_\delta,
    \end{equation}
    where
    \begin{equation*}
        \mathbf{C}_1 := \begin{bmatrix}
            \diag\{ 1, 0, 1 \} & 0 \\ 
            0 & \diag \{ 1, 0, 1 \}
        \end{bmatrix}, \quad
        \mathbf{C}_\delta := \begin{bmatrix}
            \diag\{ 1, 0, -1 \} & 0 \\
            0 & \diag\{ 1, 0, -1 \}
        \end{bmatrix}.
    \end{equation*}
    Set
    \begin{equation} \label{eq-definition-lambda-mu}
        \lambda := \frac{1}{2} \left( \coth \frac{\theta_1}{2} + \coth \frac{\theta_3}{2} \right), \quad \mu := \frac{1}{2} \left( \coth \frac{\theta_1}{2} - \coth \frac{\theta_3}{2} \right),
    \end{equation}
    so that $\mathbf C=\lambda\mathbf C_1+\mu\mathbf C_\delta$.

    To solve $\mathbf Z^*\mathbf S+\mathbf S\mathbf Z=-\mathbf C$, let $\mathbf S_1$ and $\mathbf S_\delta$ solve the Lyapunov equations with right-hand sides $-\mathbf C_1$ and $-\mathbf C_\delta$, respectively. Linearity then gives
    \begin{equation} \label{eq-boson-chain-model-full-covariance-matrix}
        \mathbf{S} = \lambda \mathbf{S}_1 + \mu \mathbf{S}_\delta.
    \end{equation}

    Stability of $\mathbf Z$ guarantees uniqueness. Direct substitution gives
    \begin{equation*}
        \mathbf{S}_1 = \1_6,
    \end{equation*}
    and 
    \begin{equation*}
        \mathbf{S}_\delta = \frac{2}{4 \omega^2 + 1} \begin{bmatrix}
            -\frac{1}{2} \diag\{ -1, 0, 1 \} & \begin{bmatrix}
                0 & -\omega & 0 \\
                \omega & 0 & -\omega \\
                0 & \omega & 0
            \end{bmatrix} \\
            \begin{bmatrix}
                0 & \omega & 0 \\
                -\omega & 0 & \omega \\
                0 & -\omega & 0
            \end{bmatrix} & -\frac{1}{2} \diag\{ -1, 0, 1 \}
        \end{bmatrix}.
    \end{equation*}
    Hence \eqref{eq-boson-chain-model-full-covariance-matrix} gives the covariance matrix $\mathbf S$.

    \begin{remark} \label{remark-symplectic-eigenvalues-by-iJS}
        The symplectic eigenvalues of a positive covariance matrix $\mathbf S$ are the positive eigenvalues of $\mi\mathbf J\mathbf S$. Indeed, \eqref{eq-symplectic-diagonalization-covariance-matrix} gives
        \begin{align}
            \mi \mathbf{J} \mathbf{S} &= \mi \mathbf{J} \mathbf{G}^* \mathbf{D}_{\coth} \mathbf{G} = \mi \mathbf{G}^{-1} \mathbf{J} \mathbf{D}_{\coth} \mathbf{G} \nonumber \\
            &= \mathbf{G}^{-1} \cdot \mi \begin{bmatrix} 
                0 & D_{\coth} \\
                - D_{\coth} & 0
            \end{bmatrix} \cdot \mathbf{G}. \label{eq-eigenvalues-iJS}
        \end{align}
        Hence the eigenvalues of $\mi\mathbf J\mathbf S$ are
        \begin{equation*}
            \pm\coth(\beta_1/2),\ldots,\pm\coth(\beta_d/2).
        \end{equation*}
    \end{remark}

    For the present covariance matrix,
    \begin{equation} \label{eq-matrix-JS-boson-chain}
        \mathbf{J} \mathbf{S} = \begin{bmatrix}
            \frac{2 \mu}{4 \omega^2 + 1} \begin{bmatrix}
                0 & \omega & 0 \\
                - \omega & 0 & \omega \\ 
                0 & - \omega & 0
            \end{bmatrix} & \diag \left\{ \lambda + \frac{\mu}{4 \omega^2 + 1}, \lambda, \lambda - \frac{\mu}{4 \omega^2 + 1}  \right\} \\ 
            - \diag \left\{ \lambda + \frac{\mu}{4 \omega^2 + 1}, \lambda, \lambda - \frac{\mu}{4 \omega^2 + 1}  \right\} & \frac{2 \mu}{4 \omega^2 + 1} \begin{bmatrix}
                0 & \omega & 0 \\
                - \omega & 0 & \omega \\ 
                0 & - \omega & 0
            \end{bmatrix}
        \end{bmatrix}.
    \end{equation}

    \begin{remark} \label{remark-special-block-matrix-unitary-transform}
        A block matrix of the form
        \begin{equation*} 
            A = \begin{bmatrix}
                X & Y \\ 
                - Y & X
            \end{bmatrix},
        \end{equation*}
        with $X,Y\in M_d(\C)$, is unitarily block diagonalized by
        \begin{equation*}
            \mathcal V := \frac{1}{2} \begin{bmatrix} 
                (1 + \mi) \1_d & (1 - \mi) \1_d \\ 
                (1 - \mi) \1_d & (1 + \mi) \1_d
            \end{bmatrix},
        \end{equation*}
        namely,
        \begin{equation} \label{eq-boson-chain-model-diagonalized-A} 
            \mathcal V^* A \, \mathcal V = \begin{bmatrix}
                X - \mi Y & 0 \\ 
                0 & X + \mi Y
            \end{bmatrix}.
        \end{equation}
        If $X$ and $Y$ are real, the two diagonal blocks are complex conjugates.
    \end{remark}

    By Remark \ref{remark-special-block-matrix-unitary-transform}, it suffices to compute the eigenvalues of
    \begin{equation} \label{eq-X-iY-from-JS}
        \begin{bmatrix}
            \mi \left( \lambda + \frac{\mu}{4 \omega^2 + 1} \right) & \frac{2 \omega \mu}{4 \omega^2 + 1} & 0 \\
            -\frac{2 \omega \mu}{4 \omega^2 + 1} & \mi \lambda & \frac{2 \omega \mu}{4 \omega^2 + 1} \\ 
            0 & -\frac{2 \omega \mu}{4 \omega^2 + 1} & \mi \left( \lambda - \frac{\mu}{4 \omega^2 + 1} \right)
        \end{bmatrix}.
    \end{equation}
    Its eigenvalues are
    \begin{equation} \label{eq-symplectic-eigenvalues-boson-chain-model}
        \mi \lambda, \quad \mi \lambda + \mi \mu \frac{\sqrt{8 \omega^2 + 1}}{4 \omega^2 + 1}, \quad \mi \lambda - \mi \mu \frac{\sqrt{8 \omega^2 + 1}}{4 \omega^2 + 1}.
    \end{equation}
    Hence the symplectic eigenvalues of $\mathbf S$ are the absolute values of the coefficients of $\mi$ in \eqref{eq-symplectic-eigenvalues-boson-chain-model}.
    
    Define
    \begin{equation*}
        r := \frac{\sqrt{8 \omega^2 + 1}}{4 \omega^2 + 1}.
    \end{equation*}
    Since $\omega\neq0$, one has $0<r<1$. Substituting \eqref{eq-definition-lambda-mu}, the three symplectic eigenvalues are
    \begin{align}
        \nu_0&:=\frac{1}{2}\coth(\theta_1/2)+\frac{1}{2}\coth(\theta_3/2), \nonumber\\
        \nu_\pm&:=\frac{1\pm r}{2}\coth(\theta_1/2)
        +\frac{1\mp r}{2}\coth(\theta_3/2). \label{eq-boson-chain-symplectic-eigenvalues}
    \end{align}
    Each is a strict convex combination of numbers greater than one, and hence $\nu_0,\nu_+,\nu_->1$. Since $\theta_1\neq\theta_3$ and $r>0$, these three values are pairwise distinct. Thus there are unique, pairwise distinct inverse temperatures $ \beta_0, \beta_+, \beta_- > 0$ such that
    \[
        \coth(\beta_\sigma/2)=\nu_\sigma,
        \qquad \sigma\in\{0,+,-\}.
    \]
    In Williamson coordinates, the inequalities $\nu_0,\nu_+,\nu_->1$ are equivalent to $\mathbf S+\mi\mathbf J>0$. Proposition \ref{proposition-unique-invariant-faithful-state} therefore shows that the semigroup has a unique faithful normal invariant state. Its Williamson decomposition has three singleton temperature blocks, so Corollary \ref{corollary-distinct-temperatures} yields a positive KMS spectral gap. Thus, this example has a positive KMS spectral gap but no GNS spectral gap.

    The standardized coefficient matrices can also be computed explicitly. One may first construct a symplectic diagonalizer $\mathbf G$ of $\mathbf S$ by the method of \cite[Appendix B]{pirandola2009correlation}, set $\mathbf M=\mathbf G^{-1}$, recover the complex-linear and conjugate-linear parts $M_1,M_2$ from \eqref{eq-reconstruction-S-from-matrix-identification}, and then apply Proposition \ref{proposition-parameters-standardized-GQMS}.
    
    \section*{Acknowledgements}
    
    The author thanks Professor Franco Fagnola and Dr. Huayu Liu for helpful discussions. The author also thanks the anonymous referees for their careful comments and for identifying errors in an earlier version of the paper.



    
    \appendix

    \section{Matrix Identification of Real-linear Operators} \label{section-appendix-matrix-identification}

    Let $S:\C^d\to\C^d$ be real-linear. There are unique complex matrices $S_1,S_2$ such that
    \begin{equation*}
        S z = S_1 z + S_2 \overline{z}, \quad \forall z \in \C^d.
    \end{equation*}
    Identifying the real Hilbert spaces $\C^d$ and $\R^{2d}$ by $z\mapsto(\Re z,\Im z)^T$, define
    \begin{equation} \label{eq-identification-real-linear-operators}
        \mathbf{S} := \begin{bmatrix}
            \Re S_1 + \Re S_2 & \Im S_2 - \Im S_1 \\ 
            \Im S_1 + \Im S_2 & \Re S_1 - \Re S_2
        \end{bmatrix}, \quad \mathbf{z} := \begin{bmatrix}
            \Re z \\ \Im z
        \end{bmatrix}.
    \end{equation}

    We call $\mathbf S$ the \emph{real matrix representation} of $S$.

    Conversely, if a real $2d\times2d$ matrix is written in $d\times d$ blocks as
    \begin{equation*}
        \mathbf{S} = \begin{bmatrix}
            S_{11} & S_{12} \\ 
            S_{21} & S_{22}
        \end{bmatrix},
    \end{equation*}
    then its associated real-linear operator is
    \begin{equation} \label{eq-reconstruction-S-from-matrix-identification}
        S z = \left( \frac{S_{11} + S_{22}}{2} + \mi \, \frac{S_{21} - S_{12}}{2}  \right) z + \left( \frac{S_{11} - S_{22}}{2} + \mi \, \frac{S_{12} + S_{21}}{2} \right) \overline{z}.
    \end{equation}

    The adjoint of $S$ with respect to $\Re\langle\cdot,\cdot\rangle$ is denoted by $S^\sharp$. It is given by
    \begin{equation} \label{eq-adjoint-real-linear-operator}
        S^\sharp z = S_1^* z + S_2^T \overline{z}, \quad \forall z \in \C^d,
    \end{equation}
    and is represented by $\mathbf S^*$. Since $\mathbf S$ is real, this is its transpose.

    The representation preserves the action and the real inner product.

    \begin{lemma}
        Let $S$ be real-linear and let $y,z\in\C^d$. Then
        \begin{equation*}
            \mathbf{S} \mathbf{z} = \begin{bmatrix}
                \Re S z \\ 
                \Im S z
            \end{bmatrix}, \quad \langle \mathbf{y}, \mathbf{z} \rangle = \Re \langle y, z \rangle.
        \end{equation*}
    \end{lemma}
    \begin{proof}
        Both identities follow directly from \eqref{eq-identification-real-linear-operators}.
    \end{proof}

    Consequently,
    \begin{equation*}
        \mathbf{S}^n \mathbf{z} = \begin{bmatrix}
            \Re S^n z \\ \Im S^n z
        \end{bmatrix}, \quad \me^{\mathbf{S}} \mathbf{z} = \begin{bmatrix}
            \Re \me^S z \\ 
            \Im \me^S z
        \end{bmatrix}.
    \end{equation*}

    A real-linear map $M:\C^d\to\C^d$ is \emph{symplectic} if
    \begin{equation*}
        \Im \langle z, w \rangle = \Im \langle M z, M w \rangle, \quad \forall z, w \in \C^d.
    \end{equation*}
    Its real matrix representation $\mathbf M$ is then symplectic in the usual sense. The matrices $\mathbf M^*$ and $\mathbf M^{-1}$ are also symplectic and represent $M^\sharp$ and $M^{-1}$, respectively.

    Writing $Mz=M_1z+M_2\overline z$, the symplectic condition is equivalent to
    \begin{align}
        M_1^* M_1 - M_2^T \overline{M_2} = \1_d, \quad M_2^T \overline{M_1} = M_1^* M_2, \label{eq-symplectic-transformation-condition-1} \\ 
        M_1 M_1^* - M_2 M_2^* = \1_d, \quad M_1 M_2^T = M_2 M_1^T. \label{eq-symplectic-transformation-condition-2}
    \end{align} 
    Moreover,
    \begin{equation} \label{eq-inverse-real-linear-operator}
        M^{-1} z = M_1^* z - M_2^T \overline{z}.
    \end{equation}

    \section{Standardization and Invariance of Spectral Gaps} \label{section-appendix-invariance-spectral-gap}

    Let $(\cT_t)_{t\geq0}$ satisfy the standing assumption. By \eqref{eq-gaussian-state-parthasarathy-form},
    \begin{equation} \label{eq-appendix-gaussian-state-parthasarathy-form}
        \rho = W(\omega) B(\mathbf{M}) \bigotimes_{j=1}^d (1 - \me^{-\beta_j}) \me^{-\beta_j N_j} B(\mathbf{M})^{-1} W(\omega)^{-1},
    \end{equation}
    where $\omega$ is the mean vector of $\rho$, $\mathbf M$ is a Williamson symplectic matrix, and $B(\mathbf M)$ is its Bogoliubov implementer. The matrix $\mathbf M$ need not be unique; for fixed $\mathbf M$, the implementer is unique up to a scalar of modulus one.

    Set
    \begin{equation} \label{eq-definition-rho-0}
        \rho_0 := \bigotimes_{j=1}^d (1 - \me^{-\beta_j}) \me^{-\beta_j N_j}.
    \end{equation}
    For a fixed choice of $\mathbf M$, define
    \begin{equation} \label{eq-construction-standardization-qms}
        \widetilde{\cT}_t (X) := \cU^* \cT_t ( \cU X \cU^*) \cU, \quad \forall X \in \cB (\mathsf{h}),
    \end{equation}
    where
    \begin{equation} \label{eq-unitary-transform}
        \cU := W(\omega) B(\mathbf{M})
    \end{equation}
    is unitary. Then $\rho_0=\cU^* \rho \, \cU$ is invariant under $(\widetilde{\cT}_t)_{t\geq0}$. We call $(\widetilde{\cT}_t)_{t\geq0}$ the \emph{standardized QMS}. Multiplying $B(\mathbf M)$ by a scalar of modulus one does not affect \eqref{eq-construction-standardization-qms}.

    The parameters of the standardized semigroup are as follows.

    \begin{proposition} \label{proposition-parameters-standardized-GQMS}
        The family $(\widetilde{\cT}_t)_{t\geq0}$ is a Gaussian QMS with unique faithful normal invariant state $\rho_0$. Its drift and diffusion operators are
        \[
            \widetilde Z=M^{-1}ZM,
            \qquad \widetilde C=M^\sharp CM.
        \]
        If $Mz=M_1z+M_2\overline z$, then its coefficient matrices are
        \begin{align}
            \widetilde{\Omega} &= M_1^* \Omega M_1 + M_2^T \Omega^T \overline{M_2} + M_1^* \kappa \overline{M_2} + M_2^T \kappa^* M_1, \label{eq-widetilde-omega} \\
            \widetilde{\kappa} &= M_1^* \Omega M_2 + M_2^T \Omega^T \overline{M_1} + M_1^* \kappa \overline{M_1} + M_2^T \kappa^* M_2, \label{eq-widetilde-kappa} \\
            \widetilde{U} &= U \overline{M_1} + \overline{V} M_2,
            \qquad \widetilde{V} = \overline{U} M_2 + V \overline{M_1}, \label{eq-widetilde-U-V} \\
            \widetilde{\zeta} &= M^{\sharp} (\zeta - 2 Z^{\sharp} J \omega) = 0. \label{eq-widetilde-zeta}
        \end{align}
    \end{proposition}
    \begin{proof}
        Unitary conjugation shows immediately that $(\widetilde{\cT}_t)_{t\geq0}$ is a QMS and that $\rho_0$ is its unique faithful normal invariant state. Its action on Weyl operators is
        \begin{align}
            \widetilde{\cT}_{t} ( W(z) ) &= \cU^* \cT_t ( \cU W(z) \cU^* ) \cU \nonumber \\
            &= \exp{- 2 \mi \Im \langle \omega, M z \rangle} \, \cU^* \cT_t ( W(Mz) ) \cU \nonumber \\ 
            &= \exp{- 2 \mi \Im \langle \omega, M z \rangle} \exp{ 2 \mi \Im \langle \omega, \me^{t Z} M z \rangle } \nonumber \\         
            &\quad \cdot \exp{ - \frac{1}{2} \int_0^t \Re \langle \me^{s Z} M z, C \me^{s Z} M z \rangle \md s + \mi \int_0^t \Re \langle \zeta, \me^{s Z} M z \rangle \md s} \nonumber \\
            &\quad \cdot W( M^{-1} \me^{t Z} M z), \label{eq-explicit-action-diagonalized-semigroup-on-weyl-operators}
        \end{align}
        Define
        \begin{equation} \label{eq-appendix-tilde-Z-and-tilde-C}
            \widetilde{Z} := M^{-1} Z M, \quad \widetilde{C} := M^{\sharp} C M.
        \end{equation}
        Then \eqref{eq-explicit-action-diagonalized-semigroup-on-weyl-operators} becomes
        \begin{align*}
            \widetilde{\cT}_t (W(z)) &= \exp{\mi \int_0^t \Re \langle \zeta, \me^{s Z} M z \rangle \md s} \exp{\mi \int_0^t \Im \langle 2 \omega, Z \me^{s Z} M z  \rangle \md s } \\
            &\quad \cdot \exp{-\frac{1}{2} \int_0^t \Re \langle \me^{s Z} M z, C \me^{s Z} M z \rangle \md s } W(\me^{t \widetilde{Z}} z) \\
            &= \exp{\mi \int_0^t \Re \langle M^\sharp ( \zeta + 2 Z^{\sharp} J^{\sharp} \omega ) , M^{-1} \me^{s Z} M z \rangle \md s} \\
            &\quad \cdot \exp{ -\frac{1}{2} \int_0^t \Re \langle M^{-1} \me^{s Z} M z, M^{\sharp} C M M^{-1} \me^{s Z} M z \rangle \md s} W(\me^{t \widetilde{Z}} z) \\
            &= \exp{ - \frac{1}{2} \int_0^t \Re \langle \me^{s \widetilde{Z}} z, \widetilde{C} \me^{s \widetilde{Z}} z \rangle \md s} W(\me^{t \widetilde{Z}} z).
        \end{align*}
        Here \eqref{eq-relationship-S-C-Z-omega-zeta} gives $\zeta=2Z^\sharp J\omega$, so the linear term vanishes and \eqref{eq-widetilde-zeta} holds. The displayed Weyl formula also proves Gaussianity.

        It remains to identify the coefficient matrices. Substitution of \eqref{eq-widetilde-U-V} into $\widetilde C=M^\sharp CM$ gives
        \begin{align*}
            \widetilde{C} z &= M^{\sharp} C M z \\
            &= M_1^* ( U^T \overline{U} + V^T \overline{V} ) M_1 z + M_2^T (U^* U + V^* V) \overline{M_2} z \\
            &\quad + M_2^T ( U^* \overline{V} + V^* \overline{U} ) M_1 z + M_1^* (U^T V + V^T U) \overline{M_2} z \\
            &\quad + M_2^T ( U^* U + V^* V ) \overline{M_1} \overline{z} + M_1^* ( U^T \overline{U} + V^T \overline{V} ) M_2 \overline{z} \\ 
            &\quad + M_1^* ( U^T V + V^T U ) \overline{M_1} \overline{z} + M_2^T ( U^* \overline{V} + V^* \overline{U} ) M_2 \overline{z} \\
            &= \left( \widetilde{U}^T \overline{\widetilde{U}} + \widetilde{V}^T \overline{\widetilde{V}} \right) z + \left( \widetilde{U}^T \widetilde{V} + \widetilde{V}^T \widetilde{U} \right) \overline{z}.
        \end{align*}
        Similarly, \eqref{eq-inverse-real-linear-operator} gives
        \begin{align*}
            \widetilde{Z} z &= M^{-1} Z M z \\
            &= \frac{1}{2} \left( M_1^* ( U^T \overline{U} - V^T \overline{V} ) M_1 z - M_2^T ( U^* U - V^* V ) \overline{M_2} z \right. \\
            &\quad - M_2^T (U^* \overline{V} - V^* \overline{U} ) M_1 z + M_1^* ( U^T V - V^T U ) \overline{M_2} z \\
            &\quad - M_2^T ( U^* U - V^* V ) \overline{M_1} \overline{z} + M_1^* ( U^T \overline{U} - V^T \overline{V} ) M_2 \overline{z} \\ 
            &\quad \left. + M_1^* ( U^T V - V^T U ) \overline{M_1} \overline{z} - M_2^T ( U^* \overline{V} - V^* \overline{U} ) M_2 \overline{z} \, \right) \\
            &\quad + \mi M_1^* \Omega M_1 z + \mi M_2^T \overline{\Omega} \overline{M_2} z + \mi M_1^* \Omega M_2 \overline{z} + \mi M_2^T \overline{\Omega} \overline{M_1} \overline{z} \\ 
            &\quad + \mi M_1^* \kappa \overline{M_2} z + \mi M_2^T \overline{\kappa} M_1 z + \mi M_1^* \kappa \overline{M_1} \overline{z} + \mi M_2^T \overline{\kappa} M_2 \overline{z} \\
            &= \left( ( \widetilde{U}^T \overline{\widetilde{U}} -\widetilde{V}^T \overline{\widetilde{V}}) /2 + \mi \widetilde{\Omega} \right) z + \left( (\widetilde{U}^T \widetilde{V} - \widetilde{V}^T \widetilde{U})/2 + \mi \widetilde{\kappa} \right) \overline{z}.
        \end{align*}
        The matrix $\widetilde\Omega$ is Hermitian and $\widetilde\kappa$ is symmetric. Thus the transformed drift and diffusion operators have the required form, with coefficients \eqref{eq-widetilde-omega}--\eqref{eq-widetilde-zeta}.
    \end{proof}

    We next verify that the criterion in Theorem \ref{theorem-main-2} is invariant under unitary changes of the Kraus family and under the choice of Williamson diagonalization.

    The Kraus operators in a GKLS representation are not unique. In particular, unitary transformations of the Kraus operators leave the dissipative pre-generator unchanged; see \cite[Remark after Definition 3]{agredo2022kossakowski} for Gaussian QMSs and \cite{parthasarathy1992quantum} for the more general setting of uniformly continuous QMSs. We show that the blockwise rank criterion in Theorem \ref{theorem-main-2} is invariant under this freedom.

    \begin{lemma}[Unitary changes of a Kraus family] \label{lemma-kraus-nonuniqueness}
        Let $L_1,\ldots,L_m$ and $\widehat L_1,\ldots,\widehat L_m$ be two minimal Kraus families related by
        \begin{equation*}
            \widehat L_\ell=\sum_{j=1}^m Q_{\ell j}L_j,
        \end{equation*}
        where $Q\in M_m(\C)$ is unitary. Then the dissipative part $\cL-\mi[H,\,\cdot\,]$ is unchanged, and the corresponding coefficient matrices satisfy
        \begin{equation} \label{eq-unitary-transform-U-V}
            \widehat U=QU,
            \qquad
            \widehat V=\overline QV.
        \end{equation}
        Consequently, for every equal-temperature block,
        \begin{equation} \label{eq-same-kernel-unitary-transform-U-V}
            \ker\widehat U_n=\ker U_n,
            \qquad
            \ker\widehat V_n=\ker V_n.
        \end{equation}
    \end{lemma}
    \begin{proof}
        For $x\in\dom\cL$, the dissipative part is
        \begin{equation*}
            -\frac12\sum_{\ell=1}^m
            \left(L_\ell^*L_\ell x-2L_\ell^*xL_\ell+xL_\ell^*L_\ell\right)
            =\sum_{\ell=1}^m
            \left(L_\ell^*xL_\ell-\frac12\{L_\ell^*L_\ell,x\}\right).
        \end{equation*}
        Unitarity of $Q$ gives
        \begin{equation*}
            \sum_{\ell=1}^m\widehat L_\ell^*\widehat L_\ell
            =\sum_{j,k=1}^m(Q^*Q)_{jk}L_j^*L_k
            =\sum_{j=1}^mL_j^*L_j,
        \end{equation*}
        and likewise
        $\sum_\ell\widehat L_\ell^*x\widehat L_\ell
        =\sum_jL_j^*xL_j$.
        Hence the dissipative part is unchanged. Comparing the coefficients in \eqref{eq-Kraus-operator} gives \eqref{eq-unitary-transform-U-V}. The invariant state, its inverse temperatures, and the blocks $I_n$ are therefore unchanged. Finally,
        \begin{equation*}
            \widehat U_n=QU_n,
            \qquad
            \widehat V_n=\overline QV_n,
        \end{equation*}
        and $Q,\overline Q$ are invertible, which proves \eqref{eq-same-kernel-unitary-transform-U-V}.
    \end{proof}

    We now consider the non-uniqueness of the Williamson diagonalizer. To avoid ambiguity, $\widetilde U,\widetilde V$ denote the coefficients obtained from $\mathbf G$, while $\widehat{\widetilde U},\widehat{\widetilde V}$ denote those obtained from $\widehat{\mathbf G}$.

    \begin{lemma}[Non-uniqueness of the Williamson diagonalization] \label{lemma-williamson-nonuniqueness}
        Suppose that 
        \begin{equation*}
            \mathbf{S} = \mathbf{G}^* \mathbf{D}_{\coth} \mathbf{G} = \widehat{\mathbf{G}}^* \mathbf{D}_{\coth} \widehat{\mathbf{G}},
        \end{equation*}
        where the same ordering of the symplectic eigenvalues is used in both diagonalizations. Let 
        \begin{equation*}
            \mathbf{M} = \mathbf{G}^{-1}, \quad \widehat{\mathbf{M}} = \widehat{\mathbf{G}}^{-1}.
        \end{equation*}
        Then changing $\mathbf G$ to $\widehat{\mathbf G}$ produces only a unitary change of coordinates within each equal-temperature block. More precisely, for every block $I_n$ there is a unitary matrix $K_n\in M_{|I_n|}(\C)$ such that
        \begin{equation*}
            \widehat{\widetilde U}_n=\widetilde U_nK_n^T,
            \qquad
            \widehat{\widetilde V}_n=\widetilde V_nK_n^T.
        \end{equation*}
        Consequently, $\widehat{\widetilde U}_n$ is injective if and only if $\widetilde U_n$ is injective, and the same equivalence holds for $\widehat{\widetilde V}_n$ and $\widetilde V_n$.
    \end{lemma}
    \begin{proof}
        Define $\mathbf K:=\widehat{\mathbf G}\mathbf G^{-1}$. Since both factors are symplectic, so is $\mathbf K$:
        \begin{equation} \label{eq-williamson-diag-1}
            \mathbf{K}^* \mathbf{J} \mathbf{K} = \mathbf{J}.
        \end{equation}
        The two Williamson decompositions also give
        \begin{equation} \label{eq-williamson-diag-2}
            \mathbf{K}^* \mathbf{D}_{\coth} \mathbf{K} = \mathbf{D}_{\coth}.
        \end{equation}
        Equation \eqref{eq-williamson-diag-1} yields $\mathbf J\mathbf K=(\mathbf K^*)^{-1}\mathbf J$, whereas \eqref{eq-williamson-diag-2} yields $\mathbf D_{\coth}^{-1}(\mathbf K^*)^{-1}=\mathbf K\mathbf D_{\coth}^{-1}$. Hence
        \begin{equation*}
            \mathbf{D}_{\coth}^{-1} \mathbf{J} \mathbf{K} = \mathbf{D}_{\coth}^{-1} (\mathbf{K}^*)^{-1} \mathbf{J} = \mathbf{K} \mathbf{D}_{\coth}^{-1} \mathbf{J}.
        \end{equation*}
        Thus $\mathbf K$ commutes with $\mathbf D_{\coth}^{-1}\mathbf J$. Since $\mathbf D_{\coth}$ commutes with $\mathbf J$,
        $(\mathbf D_{\coth}^{-1}\mathbf J)^2=-\mathbf D_{\coth}^{-2}$. Therefore,
        \begin{equation*}
            \mathbf{K} \mathbf{D}_{\coth}^{-2} = - \mathbf{K} \mathbf{D}_{\coth}^{-1} \mathbf{J} \mathbf{D}_{\coth}^{-1} \mathbf{J} = - \mathbf{D}_{\coth}^{-1} \mathbf{J} \mathbf{D}_{\coth}^{-1} \mathbf{K} = \mathbf{D}_{\coth}^{-2} \mathbf{K}.
        \end{equation*}
        The positive matrices $\mathbf D_{\coth}$ and $\mathbf D_{\coth}^{-2}$ have the same eigenspaces. It follows that
        \begin{equation} \label{eq-williamson-diag-3}
            \mathbf{K} \mathbf{D}_{\coth} = \mathbf{D}_{\coth} \mathbf{K}.
        \end{equation}

        Combining \eqref{eq-williamson-diag-2} and \eqref{eq-williamson-diag-3},
        \begin{equation*}
            \mathbf{K}^* \mathbf{K} \mathbf{D}_{\coth} = \mathbf{K}^* \mathbf{D}_{\coth} \mathbf{K} = \mathbf{D}_{\coth}. 
        \end{equation*}
        Since $\mathbf D_{\coth}$ is invertible, $\mathbf K^*\mathbf K=\1_{2d}$, so $\mathbf K$ is orthogonal. Together with symplecticity, this gives
        \begin{equation*}
            \mathbf{J} \mathbf{K} = ( \mathbf{K}^{*} )^{-1} \mathbf{J} = ( \mathbf{K} ^{-1})^{*} \mathbf{J} = \mathbf{K} \mathbf{J}.
        \end{equation*}
        Hence the real-linear map $K$ represented by $\mathbf K$ is complex-linear and unitary. Equation \eqref{eq-williamson-diag-3} shows that $K$ commutes with $D_{\coth}$. Its eigenspaces are the equal-temperature subspaces, so
        \begin{equation} \label{eq-williamson-diag-4}
            K = K_1 \oplus \cdots \oplus K_l,
        \end{equation}
        where $K_n\in M_{|I_n|}(\C)$ is unitary.

        Since $\widehat G=KG$,
        \begin{equation*}
            \widehat{M} = \widehat{G}^{-1} = G^{-1} K^{-1} = M K^*.
        \end{equation*}
        If $Mz=M_1z+M_2\overline z$, then $\widehat M_1=M_1K^*$ and $\widehat M_2=M_2K^T$. Formula \eqref{eq-widetilde-U-V} therefore gives
        \begin{equation*}
            \widehat{\widetilde U}=\widetilde U K^T,
            \qquad
            \widehat{\widetilde V}=\widetilde V K^T,
        \end{equation*} 
        and \eqref{eq-williamson-diag-4} yields
        \begin{equation*}
            \widehat{\widetilde U}_n=\widetilde U_nK_n^T,
            \qquad
            \widehat{\widetilde V}_n=\widetilde V_nK_n^T.
        \end{equation*}
        Since $K_n^T$ is invertible, right multiplication by $K_n^T$ preserves injectivity.
    \end{proof}

    Thus the criterion in Theorem \ref{theorem-main-2} is independent of the Williamson diagonalization. The actual kernel subspace may change by the invertible coordinate map $(K_n^T)^{-1}$, but its triviality does not.

    Finally, let $(\widetilde T_t^{(s)})_{t\geq0}$ and $\widetilde L^{(s)}$ denote the induced semigroup and generator associated with $(\widetilde{\cT}_t)_{t\geq0}$. Define the unitary operator
    \begin{equation*}
        \mathscr U:\cB_2(\mathsf h)\to\cB_2(\mathsf h),
        \qquad \mathscr U(Y)=\cU^*Y\cU.
    \end{equation*}

    \begin{lemma}[Unitary equivalence of the induced semigroups] \label{lemma-unitary-induced-semigroups}
        For every $s\in[0,1]$ and $t\geq0$,
        \[
            \mathscr U T_t^{(s)}\mathscr U^*=\widetilde T_t^{(s)}.
        \]
        Consequently,
        \[
            \widetilde L^{(s)}=\mathscr U L^{(s)}\mathscr U^*,
            \qquad
            \gap\widetilde L^{(s)}=\gap L^{(s)}.
        \]
    \end{lemma}
    \begin{proof}
        For $X\in\cB(\mathsf h)$, using $\rho=\cU\rho_0\cU^*$ gives
        \begin{align*}
            &\mathscr U T_t^{(s)}\mathscr U^*
            \bigl(\rho_0^{s/2}X\rho_0^{(1-s)/2}\bigr) \\
            &\quad=\cU^*T_t^{(s)}
            \bigl(\rho^{s/2}\cU X\cU^*\rho^{(1-s)/2}\bigr)\cU \\
            &\quad=\rho_0^{s/2}\cU^*\cT_t(\cU X\cU^*)\cU
            \rho_0^{(1-s)/2} \\
            &\quad=\widetilde T_t^{(s)}
            \bigl(\rho_0^{s/2}X\rho_0^{(1-s)/2}\bigr).
        \end{align*}
        The subspace $\rho_0^{s/2}\cB(\mathsf h)\rho_0^{(1-s)/2}$ is dense in $\cB_2(\mathsf h)$, so the semigroup identity holds everywhere. The generator identity and equality of the numerical gaps follow from unitary equivalence. See also \cite[Proposition C.3]{fagnola_li_2025spectral}.
    \end{proof}

    \begingroup
    \raggedright
    \setlength{\emergencystretch}{3em}
    \biburlbigskip=0mu\relax
    \printbibliography
    \endgroup

\end{document}